\documentclass[12pt,a4paper,draft]{article}
\usepackage[T2A]{fontenc}
\usepackage[english]{babel}
\usepackage[all]{xy}
\usepackage{latexsym,amsfonts,amssymb,amsmath,longtable,
mathrsfs,theorem,cite
}
\newcommand{\qed}{\hfill{$\Box$}}

\newtheorem{Lemma}{\bfseries Lemma}
\newtheorem{Theo}[Lemma]{\bfseries Theorem}
\newtheorem{Cor}[Lemma]{\bfseries Corollary}

\theoremstyle{remark}

\newtheorem{Ex}{\bfseries Example}
\newtheorem{Conj}{\bfseries Conjecture}

\title{\vspace{-1cm} \hfill{\normalsize 20D20}{
\fontfamily{cmr} \fontseries{bx} \selectfont \\ \vspace{1cm} On the pronormality of Hall subgroups}
\thanks{The work is supported by RFBR, projects 11--01--00456, 12--01--33102.}}
\date{}
\author{\bf   Evgeny. P. Vdovin, Danila. O. Revin}

\begin{document}

\maketitle

\begin{abstract}
Fix a set of primes $\pi$. A finite group is said to satisfy $C_\pi$ or, in other words, to be a
$C_\pi$-group,  if it possesses exactly one class of conjugate  $\pi$-Hall subgroups. The pronormality of $\pi$-Hall
subgroups in~$C_\pi$-groups is proven, or, equivalently, we prove that $C_\pi$ is inherited by
overgroups of $\pi$-Hall subgroups. Thus an affirmative solution to Problem 17.44(a) from the ``Kourovka notebook'' is
obtained. We also provide an example, showing that Hall subgroups in finite groups are not pronormal in general.
\end{abstract}

\section*{Introduction}

We use the term ``group'' in the meaning ``finite group''. The notation
mod CFSG means that the result is proven by using the classification of finite simple groups.

Throughout a set of primes is denoted by $\pi$, while its complement is denoted by~$\pi'$.

A subgroup $H$ of $G$ is called a {\it $\pi$-Hall subgroup},  if it is a $\pi$-group
(i.\,e. all its prime divisors are in $\pi$), while its index is not divisible by primes from
$\pi$. The notion of $\pi$-Hall subgroup generalizes the notion of Sylow $p$-subgroup and is equal to the
second notions if  $\pi=\{p\}$. The set of $\pi$-Hall subgroup of $G$ is denoted by $\operatorname{Hall}_\pi(G)$.
A subgroup is said to be a {\em Hall subgroup}, if it is a $\pi$-Hall subgroup for a set of primes $\pi$.

According to [1] we say that $G$  {\it satisfies $($belongs to the class$)$} $E_\pi$, if  $G$ possesses a
$\pi$-Hall subgroup. If $G\in E_\pi$ and every two its $\pi$-Hall subgroups are conjugate, then we say that
$G$ { \it satisfies $C_\pi$} (and we write ${G\in C_\pi}$). If ${G\in C_\pi}$ and each $\pi$-subgroup of $G$ is
included in a $\pi$-Hall subgroup, the we say that $G$ { \it satisfies} $D_\pi$ (and we write~${G\in D_\pi}$).

Properties $E_\pi$, $C_\pi$, and $D_\pi$ generalize known properties of Sylow subgroups in case of
$\pi$-Hall subgroups, but, in contrast with the Sylow properties, arbitrary group may fail to satisfy  $E_\pi$,
$C_\pi$ or $D_\pi$. Groups satisfying these properties are called $E_\pi$-, $C_\pi$- and {\it $D_\pi$-groups}
respectively.

In the theory of properties $E_\pi$, $C_\pi$, and $D_\pi$ problems of the inheriting of these properties by subgroups,
homomorphic images, and extensions are very important. These properties are studied in [1--19] (for details see
[2,\,3]). In particular, $E_\pi$ and (mod CFSG) $D_\pi$ are known to inherit by normal subgroups, while $C_\pi$ does
not in general. However, even $E_\pi$ and $D_\pi$ are not inherited by arbitrary subgroup. Consider the following
example.

\begin{Ex}
According to [4, Theorem 3], $G=\operatorname{SL}_2(16)\simeq A_1(16)$ satisfies $D_\pi$ with
$\pi=\{3,5\}$, and a subgroup $M=\operatorname{SL}_2(4)\simeq \operatorname{Alt}_5$ of order $60=2^2\cdot 3\cdot 5$ is
included into $G$ in the natural.  This subgroup does not even satisfy $E_\pi$, since it does not contain
elements (and, therefore, subgroups) of order~$15$.
\end{Ex}

A natural problem arises: what subgroups, except normal subgroups, inherit properties $E_\pi$, $C_\pi$ and~$D_\pi$?

Clearly, a $\pi$-Hall subgroup of an $E_\pi$-group is a $\pi$-Hall subgroup in each subgroup containing it, i.~e.
{\sl $E_\pi$ is inherited by overgroups of $\pi$-Hall subgroup.}

We formulate the same statements for $C_\pi$ and $D_\pi$ as conjectures.

\begin{Conj} {\rm[20, Problem 17.44(a); 21, Problem 2; 5, Conjecture 3]}
If $G\in C_\pi$ and $H\in\operatorname{Hall}_\pi(G)$,  then $M\in C_\pi$ for every subgroup $M$ such that
$H\le M\le G$.
\end{Conj}

\begin{Conj} {\rm [20, Problem 17.44(b); 21, Problem 3]}
If $G\in D_\pi$ and $H\in\operatorname{Hall}_\pi(G)$, then
$M\in D_\pi$ for every subgroup $M$ such that~${H\le M\le G}$.
\end{Conj}

The following theorem from~[5] allows to obtain a criterion for a finite group to satisfy  $C_\pi$ in terms of its
arbitrary normal series.

\begin{Theo} {\rm([5, Theorem~1] mod CFSG)}
If $G\in C_\pi$, $H\in\operatorname{Hall}_\pi(G)$, and $A\trianglelefteq G$, then~$HA\in C_\pi$.
\end{Theo}

This theorem provides a partial affirmative answer to Conjecture~1\footnote{Moreover, this theorem allows
the authors to make the conjecture for the first time.}.

Now we give an equivalent form of Conjecture~1.

According to the definition of P.Hall, a subgroup $H$ of $G$ is called {\it pronormal},  if, for every
$g\in G$, subgroups  $H$ and $H^g$ are conjugate in $\langle H, H^g\rangle$.  Classical examples of pronormal subgroups
are:

$\bullet$ normal subgroups;

$\bullet$ maximal subgroups;

$\bullet$ Sylow subgroups;

$\bullet$ Hall subgroups of solvable groups.

It is easy to see that Conjecture 1 is equivalent to the following.

\begin{Conj}
Hall subgroup of~$C_\pi$-groups are pronormal.
\end{Conj}

We can consider a stronger statement.

\begin{Conj}
Hall subgroups of every group are pronormal.
\end{Conj}

In~[22] Conjecture 4 is proven (mod CFSG) in a particular case, namely, the following conjecture is confirmed

\begin{Conj} {\rm [20, Problem 17.45(a)]}
Hall subgroups of a finite simple group are pronormal.
\end{Conj}

Conjecture 2 also can be reformulated in the spirit of Conjecture  3, if one introduces the notion of strongly
pronormal subgroup.

A subgroup $H$ of $G$ is called {\it strongly pronormal},  if, for each $K\le H$ and every $g\in G$,
$K^g$ is conjugate with a subgroup of $H$ (but not necessary with $K$) by an element from $\langle H, K^g\rangle$.
Clearly, every strongly pronormal subgroup is pronormal. All classical examples of pronormal subgroups (normal,
maximal, Sylow subgroups, and Hall subgroups of solvable groups) appear to be examples of strongly pronormal subgroups.
The following problem arises naturally: is a pronormal subgroup always strongly pronormal? The following example
provides a negative answer to this problem

 \begin{Ex}
 Let $m$ and $n$ be natural numbers and $n/2<m< n-1$. In the symmetric group $\operatorname{Sym}_n$ the pointwise
stabilizer of an $(n-m)$-element set (a subgroup~$\operatorname{Sym}_m$) is pronormal, but is not a strongly pronormal
subgroup.
\end{Ex}

In terms of strong pronormality Conjecture 2 is equivalent to the following.

\begin{Conj}
$\pi$-Hall subgroup of~$D_\pi$-groups are strongly pronormal.
\end{Conj}

Since a finite group satisfies $D_\pi$ if and only if each its composition factor satisfies to this property [6,
Theorem~7.7 (mod CFSG)], a counter example of minimal order to equivalent Conjectures 2 and 6 should be a simple
$D_\pi$-group. So the conjectures can be derived from the following conjecture.

\begin{Conj} {\rm [20, Problem 17.45(b)]}
Hall subgroups of a finite simple group are strongly pronormal.
\end{Conj}

Finally, we can consider conjecture, that  strengthen all Conjectures~1--7.

\begin{Conj}
Hall subgroups of finite groups are strongly pronormal.
\end{Conj}

In the paper we prove, by using Theorem 1 (and therefore, by using the classification of finite simple groups), for
Conjectures 1--8, formulated above, the implication $5\Rightarrow 1$ and the equivalence
$4\Leftrightarrow 8$. Thus Conjectures 1--8 are connected to each other by the following logical diagram:

$$
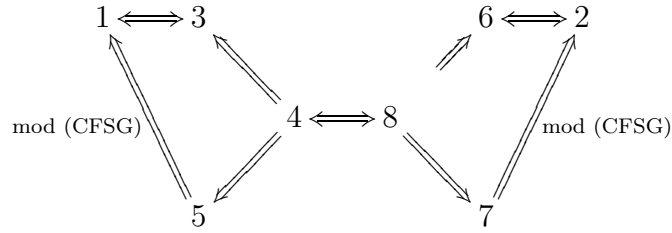
\xymatrix{
 1\ar@{<=>}[r]&3\ar@{<=}[rd]&&&6\ar@{<=}[ld
]\ar@{<=>}[r]&2\\
&&4\ar@{=>}[ld]\ar@{<=>}[r]&8\ar@{=>}[rd]&&\\
&5\ar@{=>}[luu]^{\text{mod
(CFSG)}}&&&7\ar@{=>}[ruu]_{\text{mod (CFSG)}}&}
$$

We clarify the with conjectures from the left ``wind of butterfly''. First of all, as we have already mentioned,
Conjecture 5 is true  (mod CFSG) [22, Theorem~1] and, therefore, both Conjectures 1 and 3 are true  (mod CFSG). Thus
in the paper the following theorem is proven

\begin{Theo} {\rm (mod CFSG)}
For every set of primes $\pi$ the following hold:

$(1)$~$\pi$-Hall subgroups of  $C_\pi$-groups are pronormal;

$(2)$~$C_\pi$ is inherited by overgroups of $\pi$-Hall subgroups.
\end{Theo}

Now we turn to Conjecture  4. We say that a  {\it $\pi$-conjecture} holds in $G$ for a set of primes $\pi$, if all
$\pi$-Hall subgroups of $G$ are pronormal. Thus Conjecture 4 asserts that $\pi$-conjecture holds in all finite groups
for all sets of primes~$\pi$.  $\pi$-conjecture holds in many special cases:

$\bullet$ for all simple groups [22, Theorem 1],

$\bullet$ for all groups not in $E_\pi$ (trivial),

$\bullet$ for all groups satisfying $C_\pi$ (by Theorem~2).

Thus we need to check that  $\pi$-conjecture holds in groups from
${E_\pi\setminus C_\pi}$. Notice that there exist sets $\pi$ such that
$E_\pi\setminus C_\pi=\varnothing$. Evident examples are: the set of all primes, the empty set, and an one-element set.
Every set of all primes also satisfies this property (mod CFSG) [23, Theorem A]. For such sets
$\pi$-conjecture holds in all groups. However, if there is a ``gap'' between  $E_\pi$ and $C_\pi$, then for some group
$\pi$-conjecture does not hold.

\begin{Theo}
Let a set of primes $\pi$ be such that $E_\pi\setminus C_\pi\ne\varnothing$. Then there exist $G\in E_\pi\setminus
C_\pi$ and $H\in\operatorname{Hall}_\pi(G)$ such that $H$ is not pronormal in~$G$.
\end{Theo}

It follows from the proof of Theorem 3 that we can take a regular wreath product of an arbitrary $X\in E_\pi\setminus
C_\pi$ and a cyclic group ${\Bbb Z}_p$ of order $p\in\pi'$ (here $\pi'\ne\varnothing$, since otherwise $E_\pi= C_\pi$)
as a group $G$ from the theorem.

Combining Theorems 2 and 3 we obtain

\begin{Cor} {\rm (mod CFSG)}
For every set $\pi$ of primes the following statements are equivalent:

$(1)$~$\pi$-Hall subgroups in all finite groups are pronormal;

$(2)$~$E_\pi=C_\pi$.
\end{Cor}

It would be interesting to find all sets $\pi$ such that $E_\pi=C_\pi$
(cf.~[2, Problem~7.20; 21, Problem~6]).

Theorem 3 does not just refute Conjecture 4, but also refutes a stronger Conjecture~8. Thus the implication is proven
$8\Rightarrow 4$, since this implication is made from the false statement.

Probably, the symmetry of the ``logical butterfly'' has a deeper reason: the authors do not know any counter example to
the following conjecture.

\begin{Conj}
In every group pronormal Hall subgroups are strongly pronormal.
\end{Conj}

If the conjecture is true, then each statement from the right ``wind'' will be true if and only if the symmetric
statement from the left ``wind'' is true.

Partial results on Conjectures~2 and~6 of the right ``wind'' are obtained in [24], where the conjectures are confirmed
for groups whose nonabelian factors are isomorphic to alternating, sporadic, and groups of Lie type with characteristic
in $\pi$. We discuss another possible way to solve the conjectures. As it has already mentioned, the counter example
$G$ of minimal order to any of the conjectures should be a simple  $D_\pi$-group. It is also evident that
$G$ is not a $\pi$-group and the order of $G$ is divisible by at least two primes from $\pi$
(otherwise $\pi$-Hall subgroups of $G$ are strongly pronormal). The classification of simple
$D_\pi$-groups (mod CFSG) [4, Theorem~3] implies that either $2\not\in\pi$, or $3\not\in\pi$, and
in view of  [23, Theorem B;
6, Lemma~5.1, Theorem~5.2]
$\pi$-Hall subgroups of $G$ has a Sylow tower. We recall the definition.

Let $H$ be a group and $\pi(H)=\{p_1,\dots, p_n\}$. According to
[1] we say that
$H$ {\it has a Sylow tower of complexion} $(p_1,\dots, p_n)$, if
$H$ has a normal series
$$
H=H_0>H_1>\dots>H_n=1
$$
such that each it section $H_{i-1}/H_i$ is isomorphic to a Sylow $p_i$-subgroup of~$H$.

Hall subgroups having Sylow tower of the same complexion are known to be conjugate [1, Theorem~A1]. In particular, if a
Hall subgroup has a Sylow tower, then it is pronormal.

Thus Conjecture 6 can be checked for simple $D_\pi$-groups, and, therefore, for all groups, if we can prove the
following statement.

\begin{Conj} {\rm[20, Problem 17.45(c)]}
Hall subgroups having a Sylow tower are strongly pronormal.
\end{Conj}

\section{Notations, agreements, and preliminary results}

If $H$ is a pronormal subgroup of $G$, then we write $H\operatorname{prn} G$.

\begin{Lemma}
Let $G$ be a group and $A$ be its normal subgroup.
If $G\in E_\pi$ and $H\in \operatorname{Hall}_\pi(G)$,
then $A,G/A\in E_\pi$ with ${H\cap A\in \operatorname{Hall}_\pi(A)}$,
${HA/A\in \operatorname{Hall}_\pi(G/A)}$.
\end{Lemma}

{\scshape Proof.} See [1, Lemma~1].
\qed

A finite group possessing a (sub)normal series with factors being either $\pi$- or $\pi'$-groups is called
{\it $\pi$-separable}.  Notice that a subgroup of a $\pi$-separable group is $\pi$-separable.

\begin{Lemma}
Each $\pi$-separable group satisfies~$D_\pi$.
\end{Lemma}

{\scshape Proof.} See  [7;\,1, Corollary~D5.2].
\qed

It follows from Lemma 6

\begin{Lemma}
$\pi$-Hall subgroups of a $\pi$-separable group are strongly pronormal.
\end{Lemma}

{\scshape Proof.}
Assume that $G$ is $\pi$-separable, and  $H\in\operatorname{Hall}_\pi(G)$. Let $K\le H$ and $g\in G$. Consider
$M=\langle H, K^g\rangle$. It is $\pi$-separable as a subgroup of a $\pi$-separable group, and by Lemma 6 we have
$M\in D_\pi$. It follows that a $\pi$-subgroup $K^g$ of $M$ is conjugate in $M=\langle H, K^g\rangle$ with a subgroup
of $H\in\operatorname{Hall}_\pi(M)$. Hence a subgroup  $H$ of $G$ is strongly pronormal.
\qed

We need the following statement from [5] together with Theorem 1.

\begin{Lemma} {\rm (mod CFSG)}
If $G\in C_\pi$ and  $A\trianglelefteq G$, then~${G/A\in C_\pi}$.
\end{Lemma}

{\scshape Proof.} See [5, Lemma~9].
\qed

In order to prove Theorem 2 we also use the main result from [22] that confirm Conjecture 5 and that is stated in the
following lemma.

\begin{Lemma} {\rm(mod CFSG)}
Hall subgroups of finite simple groups are pronormal.
\end{Lemma}

{\scshape Proof.} See [22, Theorem~1].
\qed

\begin{Lemma}
Let $H$ be a subgroup of $G$ and
$g\in G$,
$y\in \langle H, H^g\rangle$. Then, if $H^y$ and $H^g$ are conjugate in $\langle H^y, H^g\rangle$, then
$H$ and $H^g$ are conjugate in~$\langle H, H^g\rangle$.
\end{Lemma}

{\scshape Proof.}
Let $z\in \langle H^y, H^g\rangle$ and $H^{yz}=H^g$.  Then $z\in \langle H, H^g\rangle$, since
$\langle H^y, H^g\rangle\le\langle H, H^g\rangle $. So
$x=yz\in\langle H, H^g\rangle$ and~${H^x=H^g}$.
\qed

\begin{Lemma}
Let $\overline{\phantom{g}}:G\rightarrow G_1$ be a homomorphism,
$H\le G$.  If $H\operatorname{prn} G$, then
$\overline{H}\operatorname{prn} \overline{G}$.
\end{Lemma}

{\scshape Proof.} Clear.
\qed

\begin{Lemma}
Let $G$ be a group and $G_1,\dots, G_n$ be normal subgroups of $G$ such that $[G_i,G_j]=1$ for $i\ne
j$ and $G=G_1\cdots G_n$.  Assume that for every $i=1,\dots, n$ a pronormal subgroup $H_i$ of $G_i$ is chosen, and
$H=\langle H_1,\dots, H_n\rangle$. Then $H\operatorname{prn} G$.
\end{Lemma}

{\scshape Proof.}
Choose arbitrary element $g\in G$. Then $g=g_1\dots g_n$ for some $g_1\in G_1, \dots, g_n\in G_n$.
Since, by the condition, for every $i=1,\dots, n$, the subgroup $H_i$ is pronormal in $G_i$, there exist
$x_i\in \langle H_i,H_i^{g_i}\rangle$ such that $H_i^{x_i}=H_i^{g_i}$. In view of the condition
$[H_i,H_j]=1$ for $i\ne j$, for all $i=1,\dots ,n$ we have $H_i^g=H_i^{g_i}$. By the same arguments,
$H_i^{x_i}=H_i^{x}$, where $x=x_1\dots x_n$. It is clear that
$$
x\in \bigl\langle H_i,H_i^{g_i}\mid i=1,\dots, n\bigr\rangle
=\bigl\langle H_i,H_i^{g_i}\mid i=1,\dots,
n\bigr\rangle=\langle H, H^g\rangle.
$$
Moreover
\begin{multline}
H^g=\bigl\langle H_i^{g}\mid i=1,\dots, n\bigr\rangle
=\bigl\langle H_i^{g_i}\mid i=1,\dots, n\bigr\rangle
\\
=\bigl\langle H_i^{x_i}\mid i=1,\dots, n\bigr\rangle
\bigl\langle H_i^{x}\mid i=1,\dots, n\bigr\rangle=H^x.
\end{multline}
\qed

\begin{Lemma}
Let $G$ be a group, $H\in\operatorname{Hall}_\pi(G)$ for a set $\pi$ of primes,
$A\trianglelefteq G$ and $G=HA$. If $H\cap A\operatorname{prn} A$, then $H\operatorname{prn} G$.
\end{Lemma}

{\scshape Proof.}
Let $H\cap A\operatorname{prn} A$. Choose arbitrary $g\in G$ and show that $H^x=H^g$ for some $x\in\langle H,
H^g\rangle$. Since $G=HA$, there exist $h\in  H$ and $a\in A$ such that $g=ha$. Since $H\cap
A\operatorname{prn} A$, there exists $y\in   \langle H\cap A, H^a\cap A\rangle$ such that $H^y\cap A=H^a\cap A$. Taking
into consideration Lemma~10, in view of
$$
y\in \langle H\cap A, H^a\cap A\rangle\le\langle H, H^a\rangle
=\langle H, H^{ha}\rangle=\langle H, H^g\rangle,
$$
we may assume that $H=H^y$ and, in particular, $H\cap A=H^a\cap A=H^g\cap A$. Now $H$, $H^g$ and $g$ lie in
$N_G(H\cap A)$. Since $G=HA$, we have $G=AN_G(H\cap A)$. Notice that
$$
N_G(H\cap A)/N_G(H\cap A)\simeq AN_G(H\cap A)/A=G/A
$$
is a $\pi$-group.
Consider a normal series
$$
N_G(H\cap A)\trianglerighteq N_A(H\cap A)\trianglerighteq H\cap A\trianglerighteq 1
$$
of $N_G(H\cap A)$. Each its factor is either a  $\pi$-, or a $\pi'$-group, so
$N_G(H\cap A)$ is $\pi$-separable. By Lemma 7 we have $H\operatorname{prn} N_G(H\cap A)$. Thus
$H$ and $H^g$ are conjugate in~${\langle H,H^g\rangle}$.
\qed

\section{Proof of Theorem~2}

We prove equivalent statements (1) and (2) of the theorem simultaneously. Assume that Theorem~2 is not true and
$G\in C_\pi$ is a group of minimal order possessing a nonpronormal $\pi$-Hall subgroup $H$. Then there exists
$g\in G$ such that $M=\langle H, H^g\rangle$ does not satisfy $C_\pi$.

According to Lemma 9, $G$ is not simple. Let $A$ be a minimal normal subgroup of $G$. In view of Lemma 8 and the choice
of $G$ we have $HA/A\operatorname{prn} G/A$ and so
$H^gA=H^yA$ for some $y\in M=\langle H, H^g\rangle$. By Lemma 10 we may assume that $H=H^y$ and
$HA=H^gA$.

Now $HA\in C_\pi$, by Theorem 1, and, if $HA<G$, then the minimality of the order of $G$ implies
$H\operatorname{prn} HA$, and since $H^g\le HA$, we obtain that  $H$ and $H^g$ are conjugate in $M$
(notice that $H$ and $H^g$ are conjugate in $HA$, since
$HA\in C_\pi$).
Therefore, $G=HA$.

$A$ as a minimal normal subgroup is a direct product of isomorphic simple groups:
$$
A=S_1\times\dots\times S_n.
$$
Notice that all groups $S_i$ are nonabelian (otherwise $A$ is either a $\pi$- or a $\pi'$-group, $G=HA$ would be
$\pi$-separable, and its $\pi$-Hall subgroups are pronormal according to Lemma 7). Since all subgroups  $S_i$  are
subnormal in $G$, by Lemma~5 we have $H\cap S_i\in\operatorname{Hall}_\pi(S_i)$,
$i=1,\dots,n$. Now $H\cap S_i\operatorname{prn} S_i$ for all $i$ by Lemma 9, and
$$
H\cap A=\langle H\cap S_1,\dots, H\cap S_n\rangle\operatorname{prn} A
$$
by Lemma 12. Finally, applying Lemma 13, we conclude that $H\operatorname{prn} G$. The theorem is proven.
\qed

\section{Proof of Theorem 3}

Since $E_\pi\ne C_\pi$,  $\pi$ is not equal to the set of all primes and $\pi'\ne\varnothing$. Let $p\in\pi'$.
Assume also that $X\in E_\pi\setminus C_\pi$.  Then $X$ possesses two nonconjugate $\pi$-Hall subgroups $U$ and~$V$.

Consider the direct product $$
Y=\underbrace{X\times X\times\dots\times X\times X}\limits_{p\text{ times}}
$$
of $p$ isomorphic copies pf $X$. The map $\tau:Y\rightarrow Y$, given by
$$
(x_1,x_2,\dots,x_{p-1}, x_p)\mapsto(x_2,x_3,\dots,x_{p}, x_1),
\quad x_1,x_2,\dots,x_{p-1},x_p\in X,
$$
is an automorphism of order $p$ of $Y$. Consider the natural split extension $G$ of $Y$ by
$\langle \tau\rangle$, that is isomorphic to the regular wreath product~${X\wr\Bbb{Z}_p}$.

Since $Y$ is a normal subgroup of $G$ and the index $|G:Y|=p$ is not divisible by primes from $\pi$,  we have
$\operatorname{Hall}_\pi(G)=\operatorname{Hall}_\pi(Y)$.  In $Y$ define subgroups
$$
H=V\times \underbrace{U\times\dots\times U\times U}\limits_{p-1 \text{ times}},
\quad
K=\underbrace{U\times U\times\dots\times U}\limits_{p-1 \text{ times}}\times V
$$
in the natural way. Clearly, $H,K\in \operatorname{Hall}_\pi(Y)=\operatorname{Hall}_\pi(G)$.
Since $U$ and $V$ are not conjugate in $X$, subgroups $H$ and  $K$ are not conjugate in $Y$ and, therefore, are not
conjugate in $\langle H,K\rangle$. At the same time, by the definition of  $\tau$ we obtain the identity  $H^\tau=K$,
so $H$ and $K$ are not pronormal subgroup of  $G$. The theorem is proven. \qed

\medskip
In view of the proof of Theorem 3 note that the authors do not know any counter example to the following statement.

\begin{Conj}
Hall subgroups are pronormal in their normal closure.
\end{Conj}

\end{document}